\documentclass{article}
\usepackage{amssymb}

\usepackage{graphicx}
\usepackage{amsmath}

\newtheorem{theorem}{Theorem}

\newtheorem{corollary}[theorem]{Corollary}

\newtheorem{lemma}[theorem]{Lemma}

\newtheorem{remark}[theorem]{Remark}

\newenvironment{proof}[1][Proof]{\textbf{#1.} }{\ \rule{0.5em}{0.5em}}

\begin{document}

\title{Twist deformations in dual coordinates}
\author{Vladimir Lyakhovsky \\
Department of Theoretical Physics,\\
Sankt-Petersburg State University,\\
Ulianovskaya 1, Petrodvoretz, 198904 \\
Sankt Petersburg, Russia\\
}
\maketitle

\begin{abstract}
$\rule{7mm}{0mm}$Twist deformation $U_{%
\mathcal{F}}\left( \mathfrak{g}\right) $ is equivalent to the quantum group $%
\mathrm{Fun}_{\mathrm{def}}\left( \mathfrak{G}^{\#}\right) $ and has two
preferred bases: the one originating from $U\left( \mathfrak{g}\right) $ and
that of the coordinate functions of the dual Lie group $\mathfrak{G}^{\#}$.
The costructure of the Hopf algebra $U_{\mathcal{F}}\left( \mathfrak{g}\right)
\approx \mathrm{Fun}_{\mathrm{def}}\left( \mathfrak{G}^{\#}\right) $ is
analyzed in $\mathfrak{G}^{\#}$-group terms.
The weight diagram of the adjoint representation
of the algebra $\mathfrak{g}^{\#}$ is constructed in terms of the root
system $\Lambda \left( \mathfrak{g}\right) $.
The explicit form of the $\mathfrak{g\longrightarrow g}^{\#}$ transformation
can be obtained for any simple Lie algebra $\mathfrak{g}$ and
the factorizable chain $\mathcal{F}$ of extended Jordanian twists.
The dual group approach is used to find new solutions of the twist
equations. The parametrized family of extended Jordanian deformations $U_{%
\mathcal{EJ}}\left( \mathfrak{sl}(3)\right) $ is constructed and studied in
terms of $\mathcal{SL}\left( 3\right)^{\#}$. New realizations of the
parabolic twist are found.
\end{abstract}

\section{Introduction}

Twist deformations are usually described in terms of the initial Poincare
-Birkhoff-Witt basis borrowed from $U\left( \mathfrak{g}\right) $. This is
reasonable when the algebraic sector is of main interest.
Using the $\mathfrak{g}
$-basis we explicitly manifest that the algebraic relations in $U_{\mathcal{F%
}}\left( \mathfrak{g}\right) $ remain classical. The Lie-Poisson structure is
encoded in the twisting element $\mathcal{F}$. The result of the twist
deformation $\mathcal{F}:U\left( \mathfrak{g}\right)
\longrightarrow U_{\mathcal{%
F}}\left( \mathfrak{g}\right) $ is the Hopf algebra $U_{\mathcal{F}}
\left( \mathfrak{%
g}\right) $ with the initial multiplication, unit and counit and the
deformed universal element $\mathcal{R}_{\mathcal{F}}=\mathcal{F}_{21}%
\mathcal{F}^{-1}$ , costructure $\Delta _{\mathcal{F}}$ and antipode $S_{%
\mathcal{F}}$ \cite{Drin-1}.

Consider the Hopf algebra $U_{\mathcal{F}}\left( \mathfrak{g}\right) $ as a
quantized Lie bialgebra $\left( \mathfrak{g,g}^{\#}\right) $ , i.e. as the
deformation of $U\left( \mathfrak{g}\right) $ in
the direction of $\mathfrak{g}^{\#}$%
. According to the quantum duality principle \cite{Drin-3},\cite{Sem-Tyan}
the smooth deformation $U_{\mathcal{F}}\left( \mathfrak{g}\right) $ is the
quantum group $\mathrm{Fun}_{\mathrm{def}}\left( \mathfrak{G}^{\#}\right) $ of
the dual Lie group $\mathfrak{G}^{\#}$. The group $\mathfrak{G}^{\#}$ is the
universal covering group with the Lie algebra $\mathfrak{g}^{\#}$. For the
twisted universal enveloping
algebra $U_{\mathcal{F}}\left( \mathfrak{g}\right) $
the Lie coalgebra $\left( \mathfrak{g}^{\#}\right) ^{\ast }$ is defined by the
classical $r$-matrix
\begin{equation}
r=r^{lk}e_{l}\wedge e_{k},\qquad l,k=1,\ldots ,m  \label{r-mat}
\end{equation}
and has the composition
\begin{equation*}
\delta :\mathfrak{g\longrightarrow g\otimes g,\qquad }\delta =\mathrm{ad}%
_{\left( r\right) }\circ \Delta ^{\left( 0\right) }
\end{equation*}
Here $\Delta ^{\left( 0\right) }$ is the primitive coproduct $\Delta
^{\left( 0\right) }\left( g\right) =g\otimes 1+1\otimes g,$ for
 $g\in \mathfrak{g%
}$ and $\left\{ e_{l}|l,=1,\ldots ,m\right\} $ forms the basis of the
carrier subalgebra of the twist $\mathcal{F}$ . All this means that the
costructure of $U_{\mathcal{F}}\left( \mathfrak{g}\right) $ describes the
composition law
of $\mathrm{Fun}_{\mathrm{def}}\left( \mathfrak{G}^{\#}\right) $%
, the deformed multiplication
of $\mathrm{Fun}\left( \mathfrak{G}^{\#}\right) $.
The deformation is due to the fact that
the coordinate functions in $\mathrm{%
Fun}_{\mathrm{def}}\left( \mathfrak{G}^{\#}\right) $ do
not commute -- they are
subject to the composition rules in $U\left( \mathfrak{g}\right) $.

The structure of the dual Lie algebra $\mathfrak{g}^{\#}$ can be described in
terms of the generators of $\mathfrak{g}$ as
the ''second'' Lie structure on the
space of $\mathfrak{g}$. When the
group $\mathfrak{G}^{\#}$ is twisted by $\mathcal{F%
}$ this obviously leads not only to the deformation of the (initially
cocommutative) costructure but also to the deformation of the space of
coordinate functions. This can be explicitly seen in the Jordanian twist
deformation $\mathcal{F}_{\mathcal{J}}=e^{H\otimes \sigma }$; $\mathcal{F}_{%
\mathcal{J}}:U\left( \mathfrak{b}(2)\right) \longrightarrow U_{\mathcal{J}%
}\left( \mathfrak{b}(2)\right) $ of $U\left( \mathfrak{b}(2)\right) $ .
The algebra $%
\mathfrak{b}(2),$ with the
generators $\left\{ H,E\right\} $ and the relation $%
\left[ H,E\right] =E$, is the carrier of the twist $\mathcal{F}_{\mathcal{J}%
} $. Thus the dual algebra $\mathfrak{b}(2)^{\#}$ is
equivalent to $\mathfrak{b}(2)$%
. The classical $r$-matrix $H\wedge E$ describes the isomorphism
\begin{equation*}
r:\left\{
\begin{array}{c}
H^{\ast }\Leftrightarrow E \\
E^{\ast }\Leftrightarrow -H
\end{array}
\right\} ,
\end{equation*}
But the universal element for the Jordanian deformation $\mathcal{R}=%
\mathcal{F}_{21}\mathcal{F}^{-1}=e^{\sigma \otimes H}e^{-H\otimes \sigma }$
(with $\sigma =\ln \left( 1+E\right) $) shows that in the global dualization
of the algebras $U\left( \mathfrak{b}(2)\right) $ and $U\left( \mathfrak{b}%
(2)^{\#}\right) $ the generator $E$ is only the first term in the expansion
of $H^{\ast }$ and in fact $\sigma $ is the element dual to $H$ . This is
also clearly seen in the costructure of $U_{\mathcal{J}}\left( \mathfrak{b}%
(2)\right) =\mathrm{Fun}_{\mathrm{def}}\left( \mathfrak{B}(2)^{\#}\right) $
\begin{eqnarray}
\Delta _{\mathcal{J}}\left( H\right) &=&H\otimes e^{-\sigma }+1\otimes H,
\label{borel-group} \\
\Delta _{\mathcal{J}}\left( \sigma \right) &=&\sigma \otimes 1+1\otimes
\sigma .  \notag
\end{eqnarray}
These are the multiplication rules of the group $\mathfrak{B}(2)^{\#}$ (the
semidirect product of the ''rotation'' $\left( e^{-\sigma }\right) ^{\ast }$
and ''translation'' $H^{\ast }$ ) with the Cartan coordinate $\sigma $ that
together with $H$ forms the $\mathfrak{g}^{\#}$-basis here.

As we have seen in the example above the (dual) $\mathfrak{g}^{\#}$-coordinates
are easily obtained in the case when the carrier subalgebra of the twist $%
\mathfrak{g}_{c}$ coincides with $\mathfrak{g}$.

In most of the interesting cases the carrier is a proper subalgebra of $%
\mathfrak{g}$, $\mathfrak{g}_{c}\subset \mathfrak{g,}$ and
we have the nontrivial decomposition of the
space $V_{\mathfrak{g}}=V_{\mathfrak{g}_{c}}\oplus V_{\mathfrak{a}}$%
. In particular this is true when $\mathfrak{g}$ is a semisimple Lie algebra
\cite{Oh-Ela}. On the
subspace $V_{\mathfrak{g}_{c}}$ and in the Hopf algebra $%
U_{\mathcal{F}}\left( \mathfrak{g}_{c}\right) $ the
 $\mathfrak{g}^{\#}$-coordinates
can be obtained by using the twisting
element $\mathcal{F}$ as a map $\mathfrak{%
g\longrightarrow g}^{\#}$. The twist deformation
on the space $V_{\mathfrak{a}}$
is usually described by the set of coproducts for the basic elements of $%
\mathfrak{g}$ and are expressed in terms of a mixed basis containing elements
from both $\mathfrak{g}$- and $\mathfrak{g}_{c}^{\#}$-bases.
In this representation
it is very difficult to study the properties of the costructure.

The deformed costructure is essential not only by itself. It is necessary to
know its properties to find further deformations of $U_{\mathcal{F}}\left(
\mathfrak{g}\right) $ and to construct new twists by enlarging the element $%
\mathcal{F}$ with the additional factors, that is to find factorized
solutions of the twist equations \cite{Drin-1}. In particular the subspaces
in $V_{\mathfrak{a}}$ with quasiprimitive costructure are of great importance
\cite{Kul-Lya-Olm}, \cite{Kul-Lya}.

In this paper we shall consider universal enveloping algebras $U_{\mathcal{F}%
}\left( \mathfrak{g}\right) $ for semisimple
Lie algebra $\mathfrak{g}$ twisted by
chains of extended twists -- the factorized 2-cocycle $\mathcal{F}=\mathcal{F%
}_{p}\mathcal{F}_{p-1}\ldots \mathcal{F}_{1}$ where each $\mathcal{F}_{q}$
is the solution of the twist equations for $U_{\mathcal{F}_{q-1}\ldots
\mathcal{F}_{1}}\left( \mathfrak{g}\right) $. These solutions (called the
twisting factors \cite{L-BTF}) are of two types: Jordanian \cite{OGI} and
extension \cite{KLM}. In chains of twists the Reshetikhin twisting factor is
used only to ''rotate'' the Jordanian twist and can be always combined with
the corresponding Jordanian factor. We study the group $\mathfrak{G}_{\left(
\mathcal{F}\right) }^{\#}$ and define the weight diagram for the adjoint
representation of its Lie algebra $\mathfrak{g}^{\#}$. In Section 4 we develop
the algorithm for constructing the $\mathfrak{g}^{\#}$-coordinates in the Hopf
algebra $U_{\mathcal{F}}\left( \mathfrak{g}\right) $. In Section 5 the dual
group approach is applied to study the properties of the Hopf algebra $%
U_{\mathcal{F}}\left( \mathfrak{sl}\left( 3\right) \right) $. It is
demonstrated that the dual representation of the costructure can be used to
find new twist deformations.

\section{Dual Lie algebra}

Consider the simple Lie algebra $\mathfrak{g}$ and the Cartan-Weil basis $%
\left\{ e_{i}\mid i=1,\ldots ,n\right\} $ correlated with the decomposition $%
V_{\mathfrak{g}}=V_{\mathfrak{h}}
\oplus V_{\mathfrak{n}_{+}}\oplus V_{\mathfrak{n}_{-}}=V_{%
\mathfrak{h}}\oplus _{\left( \nu \in \Lambda ^{+}\right) }
\left( V_{\nu }\oplus
V_{-\nu }\right) $,
\begin{eqnarray*}
V_{\mathfrak{n}_{+}}\oplus V_{\mathfrak{n}_{-}}
&=&\oplus _{\left( \nu \in \Lambda
^{+}\right) }\left( V_{\nu }\oplus V_{-\nu }\right) = \\
&=&V_{\mu }\oplus _{\left( \nu \in \Lambda ^{+}\mid \nu \left( h_{\mu
}\right) =1/2\right) }\left( V_{\nu }\oplus V_{-\nu }\right) \oplus _{\left(
\xi \in \Lambda ^{+}\mid \xi \left( h_{\mu }\right) =0\right) } V_{\xi
} \oplus V_{-\mu }= \\
&=&V_{1}\oplus V_{1/2}\oplus V_{0}\oplus V_{-1/2}\oplus V_{-1}.
\end{eqnarray*}

The carrier subalgebra $\mathfrak{g}_{c}$ of the chain $\mathcal{F}$ is the
proper subalgebra of $\mathfrak{g}$, $\mathfrak{g}_{c}\subset \mathfrak{g}$.
The
generators $\left\{ e_{l}\mid l=1,\ldots ,m\right\} $ of $\mathfrak{g}_{c}$
contain the root vectors $\left\{ e_{\phi }|\phi \in \Lambda _{c}\right. $ $%
\left. \subset \Lambda \left( \frak{g}\right) \right\} $
 for the fixed subset of roots $%
\Lambda _{c}$ \cite{AKL}. Consider the subspace $V_{\mathfrak{a}}\subset V_{%
\mathfrak{g}}$ generated by $\left\{ h_{\left\{ \lambda _{0}\right\} }^{\perp
},e_{\nu }|\nu \in \Lambda \setminus \Lambda _{c}\right\} $ where the Cartan
elements $h_{\left\{ \lambda _{0}\right\} }^{\perp }$ are orthogonal to any
initial root $\lambda _{0}$ \cite{Kul-Lya-Olm} in $\mathcal{F}$. Then we
have the direct sum
decomposition $V_{\mathfrak{g}}=V_{\mathfrak{g}_{c}}\oplus V_{%
\mathfrak{a}}$. Construct the
basis $\left\{ e_{i}^{\ast }\right\} $ canonically
dual to $\left\{ e_{i}\right\} $, i.e. $\left( e_{i}^{\ast },e_{j}\right)
=\delta _{ij}$. Let $\Gamma $ be the weight diagram of the adjoint
representation $\mathrm{ad}_{\mathfrak{g}}$. Our aim is to define the weight
diagram $\Gamma ^{\#}$ for the adjoint
representation of $\mathrm{ad}_{\mathfrak{%
g}^{\#}}$. The structure of $\Gamma ^{\#}$ is described in the following
statement:

\begin{lemma}
Let $\mathfrak{g}$ be the simple Lie algebra with the root system $\Lambda
\subset \Gamma _{\mathrm{ad}\left( \mathfrak{g}\right) }$ and the Cartan
decomposition $V_{\mathfrak{g}}=
V_{\mathfrak{h}}\oplus _{\left( \nu \in \Lambda
^{+}\right) }\left( V_{\nu }\oplus V_{-\nu }\right) $ , $\mathcal{F}=%
\mathcal{F}_{p}\mathcal{F}_{p-1}\ldots \mathcal{F}_{1}$ -- the solution of
the twist equations where each $\mathcal{F}_{q}$ is a Jordanian or an
extension twisting factor, $r_{\mathcal{F}}=\sum_{l,k=1,\ldots ,\dim \left(
\mathfrak{g}_{c}\right) }r^{lk}e_{l}
\wedge e_{k}\in \mathfrak{g}\wedge \mathfrak{g}$ --
the corresponding skew solution of the
CYBE, $\left( \mathfrak{g},\mathfrak{g}%
^{\#}\right) $ -- the Lie bialgebra. Then the weight diagram $\Gamma ^{\#}$
of the dual Lie algebra $\mathfrak{g}^{\#}$ is
the union $\Gamma ^{\#}:=\Gamma _{%
\mathrm{ad}}\left( \mathfrak{g}^{\#}\right) =\Gamma _{c}^{\#}\cup \Gamma
_{a}^{\#}$ of the weight diagram $\Gamma _{c}^{\#}:=\Gamma _{\mathrm{ad}%
}\left( \mathfrak{g}_{c}^{\#}\right) $ of the carrier subalgebra $\mathfrak{g}%
_{c}^{\#}\left( \mathcal{F}\right) $ and the set of weights $\Gamma
_{a}^{\#}:=\Gamma \left( \mathfrak{g}_{a}^{\#}\right) $ . The system $\Gamma
_{c}^{\#}$ is the set of weights $\left\{ -\gamma _{i}\right\} $ opposite to
the weights $\left\{ \gamma _{l}|l=1,\ldots ,
\dim \left( \mathfrak{g}_{c}\right)
\right\} $ of the carrier $\mathfrak{g}_{c}$. For each component $r_{\phi \psi
}\neq 0$ in $r$ the generator $\left( e_{\varphi }\right) ^{\ast }$ takes
the weight $-\psi $ and $\left( e_{\psi }\right) ^{\ast }$ -- the weight $%
-\varphi $. The diagram $\Gamma _{a}^{\#}$ is the subset of the initial
weight diagram: $\Gamma _{a}^{\#}=\Gamma _{\mathrm{ad}}\left( \mathfrak{g}%
\right) \diagdown \Gamma _{\mathrm{ad}}\left( \mathfrak{g}_{c}\right) $. The
weight vectors $\left\{ e_{\eta }^{\ast }|\eta \in \Gamma _{a}^{\#}\right\} $
generate the Abelian ideal $\mathfrak{g}_{a}^{\#}$ and the dual algebra is the
semidirect sum $\mathfrak{g}^{\#}
\approx \mathfrak{g}_{c}^{\#}\vdash \mathfrak{g}_{a}^{\#}$.
\end{lemma}

\begin{proof}
Consider $\delta \left( g\right) $ for an arbitrary $g\in \mathfrak{g}:$%
\begin{eqnarray}
\delta \left( g\right) &=&\left[ r^{kl}e_{k}\wedge e_{l},g\otimes 1+1\otimes
g\right] =  \notag \\
&=&r^{kl}\left[ e_{k}g\right] \wedge e_{l}
+r^{kl}e_{k}\wedge \left[ e_{l}g\right].  \label{coalg-1}
\end{eqnarray}
Any nontrivial coproduct $\delta \left( g\right) $ has a comultiplier
belonging to $\mathfrak{g}_{c}$.
If $g\in \mathfrak{g}_{c}$ then $%
\delta \left( g\right) $ belongs
to $\mathfrak{g}_{c}\otimes \mathfrak{g}_{c}$
because $\mathfrak{g}_{c}$ is an algebra.
This proves that $\mathfrak{g}_{a}^{\#}$ is the Abelian ideal.

The classical $r$-matrix performs the deformation of the Lie bialgebra $%
r:\left( \mathfrak{g},\mathrm{Ab}\right) $
$\longrightarrow \left( \mathfrak{g},%
\mathfrak{g}^{\#}\right) $ \cite{K-L}. Each term in $r$ gives rise to a
deforming function for the original Abelian composition of the dual algebra.
Let $r$ be the classical $r$-matrix corresponding to the chain of extended
twists $\mathcal{F}$. To find the structure of the algebra $\mathfrak{g}%
^{\#} $ it is sufficient to consider the deforming functions originating
from two basic twisting factors: Jordanian and extension.

\begin{itemize}
\item  The Jordanian BTF has the canonical form $\mathcal{F}_{J}=e^{h\otimes
\sigma _{\mu }}$, $\sigma _{\mu }=\ln \left( 1+e_{\mu }\right) $ \cite{OGI}.
It can be deformed by the previous twisting factors \cite{Kul-Lya},\cite{AKL}
and reveal a more complicated structure. Still in all the cases one can
write it in the form $e^{h\otimes \widetilde{\sigma }_{\mu }}$ with $%
\widetilde{\sigma }_{\mu }=\ln \left( 1+e_{\mu }+f\left( e_{\alpha },\xi
_{q}\right) \right) $ ($\xi _{q}$ are the deformation parameters of the
previous twisting factors). The corresponding $r$-matrix is
\begin{equation*}
r_{J}=h\wedge e_{\mu },
\end{equation*}
with the Cartan element $h=h_{\mu }+\gamma h_{\perp }\subset V_{\mathfrak{h}}$\
, $\mu \left( h_{\bot }\right) =0$ and $\left[ h,e_{\mu }\right] =\mu \left(
h\right) e_{\mu };\mu \left( h\right) =\mu \left( h_{\mu }\right) =1$. The
elements $h^{\ast }$ and $h_{\perp }^{\ast }$ are canonically dual: $%
\left\langle h^{\ast },h_{\perp }\right\rangle =0,\left\langle h^{\ast
},h\right\rangle =\left\langle h_{\perp }^{\ast },h_{\perp }\right\rangle =1$.
On the carrier subalgebra $\mathfrak{g}_{\mathcal{J}c}
=\mathfrak{b}\left( 2\right) $
(equivalent to its dual $\mathfrak{g}_{\mathcal{J}c}^{\#}$) we get
\begin{eqnarray*}
\delta \left( h\right)
&=&\left[ r_{J},\Delta^0 \left( h \right) \right] =-h\wedge e_{\mu }, \\
\delta \left( e_{\mu }\right)
&=&\left[ r_{J},\Delta^0 \left( e_{mu} \right)e_{\mu }\right] =0.
\end{eqnarray*}
This means that the relations in $\mathfrak{g}_{\mathcal{J}c}^{\#}$ are defined
by the deforming function
\begin{equation*}
f_{J}^{\#}\left( e_{\mu }^{\ast },h^{\ast }\right) =h^{\ast },
\end{equation*}
and $e_{\mu }^{\ast }$ is the Cartan element in $\mathfrak{g}_{\mathcal{J}%
c}^{\#} $.

Consider $a_{\nu }\in V_{\mathfrak{a}}$
\begin{eqnarray*}
\delta \left( a_{\nu }\right) &=&\left[ r_{J},
\Delta^0 \left( a_{nu} \right)\right] = \\
&=&\nu \left( h_{\mu }+\gamma h_{\perp }\right) a_{\nu }\wedge e_{\mu
}+h\wedge \left[ e_{\mu },a_{\nu }\right] =\\
&=&\left( \frac{1}{2}+\gamma \nu \left( h_{\perp }\right) \right) a_{\nu
}\wedge e_{\mu }+C_{\mu \nu }^{\rho }h\wedge e_{\rho }
\end{eqnarray*}
The corresponding deforming function is
\begin{eqnarray*}
f_{J}^{\#}\left( e_{\mu }^{\ast },a_{\nu }^{\ast }\right) &=&-\nu \left(
h_{\mu }+\gamma h_{\perp }\right) a_{\nu }^{\ast }, \\
f_{J}^{\#}\left( h^{\ast },e_{\rho }^{\ast }\right) &=&C_{\mu \nu }^{\rho
}a_{\nu }^{\ast },
\end{eqnarray*}
i.e. $h^{\ast }$ has the weight ( $-\mu $ ) with respect to the original root
system $\Lambda $.
\begin{eqnarray*}
\delta \left( a_{-\mu }\right)
&=& e_{\mu }\wedge a_{-\mu }-2\gamma h\wedge h_{\perp };
\end{eqnarray*}
\begin{eqnarray*}
f_{J}^{\#}\left( e_{\mu }^{\ast },a_{-\mu }^{\ast }\right) &=&a_{-\mu
}^{\ast }, \\
f_{J}^{\#}\left( h^{\ast },h_{\perp }^{\ast }\right) &=&-2\gamma a_{-\mu
}^{\ast },
\end{eqnarray*}

Notice that $h^{\ast }$\ does not ''shift'' the element $h_{\perp }^{\ast }$%
\ if $h$\ is proportional to $h_{\mu }$.

Summing up we
find that the adjoint action of the algebra $\mathfrak{g}_{\mathcal{J}c}^{\#}$
is fixed by the relations
\begin{equation*}
\begin{array}{ll}
f_{J}^{\#}\left( e_{\mu }^{\ast },h^{\ast }\right) =h^{\ast }, &
f_{J}^{\#}\left( h^{\ast },a_{\nu }^{\ast }\right) =C_{\mu \nu }^{\rho
}a_{\nu }^{\ast }, \\
f_{J}^{\#}\left( h^{\ast },h_{\perp }^{\ast }\right) =-2\gamma a_{-\mu
}^{\ast }, & f_{J}^{\#}\left( e_{\mu }^{\ast },a_{-\mu }^{\ast }\right)
=a_{-\mu }^{\ast }. \\
f_{J}^{\#}\left( e_{\mu }^{\ast },a_{\pm \nu }^{\ast }\right) =\mp \nu
\left( h\right) a_{\pm \nu }^{\ast }, &
\end{array}
\end{equation*}
In this case the deforming function $f_{J}^{\#}$ is itself the Lie
composition, it fixes the Lie algebra $\mathfrak{g}^{\#}$ and
$\mathfrak{g}_{\mathcal{J}c}^{\#}$\ is a
subalgebra, $\mathfrak{g}^{\#}\supset \mathfrak{g}_{%
\mathcal{J}c}^{\#}$. The adjoint action
of $\mathfrak{g}_{c}^{\#}$ on $V^{\ast }$
can be decomposed into the direct
sum: $\mathrm{ad}\left( \mathfrak{g}_{\mathcal{J%
}c}^{\#}\right) |_{V^{\ast }}=\mathrm{ad}\left( \mathfrak{g}_{\mathcal{J}%
c}^{\#}\right) |_{V_{\mathfrak{g}_{\mathcal{J}c}}^{\ast }}
\oplus d\left( \mathfrak{g}%
_{\mathcal{J}c}^{\#}\right) |_{V_{\mathfrak{a}}^{\ast }}$ (with the
representation $d\left( \mathfrak{g}_{\mathcal{J}c}^{\#}\right) $ describing
the action on $V_{\mathfrak{a}}^{\ast }$).

\item  The extension BTF. Consider the $r$-matrix
\begin{equation*}
r_{E}=e_{\mu }\wedge e_{\nu },\quad \nu \neq -\mu .
\end{equation*}
(Remember that the element itself must not satisfy the CYBE. Only the full $%
r $-matrix must be the solution.)

\begin{eqnarray*}
\delta \left( h\right) &=&\left[ r_{E},
\Delta^0 \left( h \right)\right]= \\
&=&-\left( \mu +\nu \right) \left( h\right) e_{\mu }\wedge e_{\nu }, \\
\delta \left( e_{\xi }\right) &=&\left[ r_{E},
\Delta^0 \left( e_{\xi} \right)\right] \\
&=&C_{\mu \xi }^{\rho }e_{\rho }\wedge e_{\nu }+C_{\nu \xi }^{\chi }e_{\mu
}\wedge e_{\chi }, \\
\delta \left( e_{-\mu }\right) &=&\left[ r_{E},
\Delta^0 \left( e_{-\mu } \right)\right] = \\
&=&2h_{\mu }\wedge e_{\nu }+C_{\nu ,-\mu }^{\rho }e_{\mu }\wedge e_{\rho },
\end{eqnarray*}

The deforming function $f_{E}^{\#}$ looks as follows:

\begin{eqnarray*}
f_{E}^{\#}\left( e_{\mu }^{\ast },e_{\nu +\xi }^{\ast }\right) &=&C_{\nu \xi
}^{\nu +\xi }e_{\xi }^{\ast },\qquad f_{E}^{\#}\left( e_{\mu +\xi }^{\ast
},e_{\nu }^{\ast }\right) =C_{\lambda \nu }^{\mu +\xi }e_{\xi }^{\ast }, \\
f_{E}^{\#}\left( e_{\mu }^{\ast },e_{\nu }^{\ast }\right) &=&-\left( \mu
+\nu \right) \left( h\right) h^{\ast },\qquad f_{E}^{\#}\left( e_{\mu
}^{\ast },e_{\nu -\mu }^{\ast }\right) =C_{\nu ,-\mu }^{\nu -\mu }e_{-\mu
}^{\ast }, \\
f_{E}^{\#}\left( h_{\mu }^{\ast },e_{\nu }^{\ast }\right) &=&2e_{-\mu
}^{\ast }.
\end{eqnarray*}
If we assign to the basic vector $e_{\tau }^{\ast }$ the root $\tau $ then $%
\mathrm{ad}\left( e_{\mu }^{\ast }\right) $ and $\mathrm{ad}\left( e_{\nu
}^{\ast }\right) $\ act on $e_{\tau }^{\ast }$ as having the roots $-\nu $
and $-\mu $ correspondingly. When $h$ becomes proportional to $h_{\mu }$ (or
correspondingly to $h_{\nu }$) the operator $\mathrm{ad}\left( e_{\nu
}^{\ast }\right) $ ( $\mathrm{ad}\left( e_{\lambda }^{\ast }\right) $ )
cannot shift $h_{\perp }^{\ast }$ .
\end{itemize}

Combining the properties of deforming functions corresponding to
different twisting factors one can check that for the
integral composition $f^{\#}=\sum f_{q}^{\#}$ the chain $F$
of extended twists the carrier $g_{c}^{\#}$ is a subalgebra in $%
g^{\#}$ its space is the direct sum $V_{\mathfrak{g}}^{\ast }=V_{\mathfrak{g}%
_{c}}^{\ast }\oplus V_{\mathfrak{a}}^{\ast }$ and the dual algebra is a
semidirect sum $g^{\#}\approx g_{c}^{\#}\vdash g_{a}^{\#}$.%
\end{proof}

Notice that in the general case the dual algebra $g^{\#}$\ must not contain $%
g_{c}^{\#}$\ as a subalgebra.

When the adjoint action $\mathrm{ad} \left( \mathfrak{g}_{c}^{\#}\right) $\%
is considered on $h_{\perp }^{\ast }$\ it\ always behaves as ''orthogonal''
to $h^{\ast }$.

\section{Second classical limit for twisted algebras}

In the twisted universal enveloping algebra considered as the deformed
algebra of functions $U_{\mathcal{F}}\left( \mathfrak{g}\right) =\mathrm{Fun}_{%
\mathrm{def}}\left( \mathfrak{G}^{\#}\right) $ the group $\mathfrak{G}^{\#}$ is
the universal enveloping group realized in terms of formal power series.

We have seen in the previous Section how the relations of the dual algebra $%
\mathfrak{g}^{\#}$ are encoded in the weight diagram $\Gamma ^{\#}$. This
description refers to the basis of functionals $e_{\mu}^{\ast }$ canonically
dual to the generators $\left\{ e_{\mu}\right\} $ of $\mathfrak{g}$. Thus we can
construct the group compositions for $\mathfrak{G}^{\#}$ as the coproducts in
the Hopf algebra $\mathrm{Fun}\left( \mathfrak{G}^{\#}\right) $ in terms of $%
e_{\mu}$ considered as exponential coordinates for this group. Unfortunately
only in a few trivial cases the transformation to dual (noncommuting)
coordinates can be determined by the direct comparison of $\mathrm{Fun}\left(
\mathfrak{G}^{\#}\right) $ with
 $\mathrm{Fun}_{\mathrm{def}}\left( \mathfrak{G}%
^{\#}\right) =U_{\mathcal{F}}\left( \mathfrak{g}\right) $. In the general
situation the Hopf algebra $U_{\mathcal{F}}\left( \mathfrak{g}\right) $ in its
initial form is inappropriate to extract the transformation $\mathfrak{g}$ $%
\longrightarrow \mathfrak{g}^{\#}$.

In the Hopf
algebra $\mathrm{Fun}_{\mathrm{def}}\left( \mathfrak{G}^{\#}\right) $
the composition law of the group $\mathfrak{G}^{\#}$ is deformed (by the
noncommutativity of the coordinates). In quantum deformation the
compositions are usually described in terms of undeformed coordinates, i. e.
the proper coordinates of undeformed $\mathrm{Fun}\left( \mathfrak{G}%
^{\#}\right) $ are sufficient to describe the multiplication of the deformed
$\mathrm{Fun}_{\mathrm{def}}\left( \mathfrak{G}^{\#}\right) $.
Thus we must find
the Hopf algebra $\mathrm{Fun}\left( \mathfrak{G}^{\#}\right) $ with the
commutative multiplication law. The transition $\mathrm{Fun}_{\mathrm{def}%
}\left( \mathfrak{G}^{\#}\right)
\longrightarrow \mathrm{Fun}\left( \mathfrak{G}%
^{\#}\right) $ is called the second classical limit \cite{Drin-3}. It was
proved \cite{L-dual} that for the standard quantization (with the parameter $%
q=e^{h}$) the second classical limit can be obtained by the trivial scaling
of the Lie algebra generators by $\varepsilon $ and tending $\varepsilon
,h\longrightarrow 0$\ provided that $\frac{h}{\varepsilon }=\mathrm{const}$.
Taking into account that the transition to the second classical limit is a
general procedure that does not depend on the particular form of
quantization (of a Lie bialgebra) we come to the conclusion that the same
algorithm is to be applied in the case of twisted universal enveloping
algebras $U_{\mathcal{F}}$.

We are dealing with the factorizable twists $\mathcal{F}=\mathcal{F}_{p}%
\mathcal{F}_{p-1}\ldots \mathcal{F}_{1}$ with Jordanian and extension
twisting factors $\mathcal{F}_{j}$. For the sequence $\mathcal{F}$ one can
attribute the fixed number of deformation parameters $\left\{ \xi
_{s=1,\ldots ,l}\mid l\leq p\right\} $, each of them corresponds to an
automorphism of $\mathfrak{g}$. This number is equal to the number of
independent Jordanian factors $\mathcal{F}_{\mathcal{J}s}$. This means that
there exists the parametrized solution $r\left( \left\{ \xi _{s}\right\}
\right) $ of the CYBE, and as a result
we get the parametrized sets $\mathfrak{g}%
^{\#}\left( \left\{ \xi _{s}\right\} \right) $ and $\mathfrak{G}^{\#}\left(
\left\{ \xi _{s}\right\} \right) $.

The automorphism that injects the deformation parameters together with
the mentioned above scaling leads to the substitution
\begin{equation}
e_{j}\longrightarrow \frac{\xi _{s}^{\vartheta _{\left( s\right) j}}}{%
\varepsilon }\widehat{e_{j}}.  \label{subst}
\end{equation}
(Here $\vartheta _{\left( s\right) }$ is the natural grading corresponding
to the $s$-th link of the chain $\mathcal{F}$.) The classical limit
algorithm prescribes that the parameters $\varepsilon $ and $\xi _{s}$ are
to be turned to zero while their relations are fixed
\begin{eqnarray*}
\varepsilon ,\xi _{s} &\longrightarrow &0, \\
\frac{\xi _{s}}{\varepsilon } &=&\zeta _{s}.
\end{eqnarray*}
The result is the algebra of functions over the group $\mathfrak{G}^{\#}\left(
\left\{ \zeta _{s}\right\} \right) $,
\begin{equation}
\lim_{\varepsilon ,\xi _{s}\longrightarrow 0;\quad \xi _{s}/\varepsilon
=\zeta _{s}}U_{\mathcal{F}\left( \left\{ \xi _{s}\right\} \right) }\left(
\mathfrak{g}\left( \varepsilon \right) \right)
=\mathrm{Fun}\left( \mathfrak{G}%
^{\#}\left( \left\{ \zeta _{s}\right\} \right) \right) .  \label{cl-lim}
\end{equation}

In particular one can check that when $\mathcal{F}$ is a canonically
parametrized chain of extended twists in the limit (\ref{cl-lim}) the
costructure constants of $U_{\mathcal{F}\left( \left\{ \xi _{s}\right\}
\right) }\left( \mathfrak{g}\left( \varepsilon \right) \right) $ remains
finite, the multiplication law becomes Abelian and the Lie algebra
of the obtained group coincides
with $\mathfrak{g}^{\#}\left( \left\{ \zeta _{s}\right\} \right)
.$ The corresponding calculations are performed in Appendix.

\section{Group coordinates in terms of dual algebra coordinates}

In the case of the factorizable twist $\mathcal{F}=\mathcal{F}_{p}\mathcal{F}%
_{p-1}\ldots \mathcal{F}_{1}$ the dual group $\mathfrak{G}^{\#}$ is solvable.
Thus the exponential coordinates are the most suitable to describe $\mathrm{%
Fun}_{\mathrm{def}}\left( \mathfrak{G}^{\#}\right)
=U_{\mathcal{F}}\left( \mathfrak{g%
}\right) $. In these coordinates the Hopf algebra $\left( \mathrm{Fun}_{%
\mathrm{def}}\left( \mathfrak{G}^{\#}\right) \right) ^{\ast }$ dual
to $\mathrm{%
Fun}_{\mathrm{def}}\left( \mathfrak{G}^{\#}\right) $ takes the form of the
universal enveloping algebra $U\left( \mathfrak{g}_{\left( \mathfrak{g}\right)
}^{\#}\right)$ $ =\left( \mathrm{Fun}_{\mathrm{def}}
\left( \mathfrak{G}^{\#}\right)
\right) ^{\ast }$. (Here index $\left( \mathfrak{g}\right) $ indicates that the
Hopf algebra $U\left( \mathfrak{g}_{\left( \mathfrak{g}\right) }^{\#}\right) $
is obtained in terms of (the initial) $\mathfrak{g}$-basis.)

The algebra $\mathfrak{g}^{\#}$ is fixed by the
classical $r$-matrix. One can construct the universal
enveloping $U\left( \mathfrak{g}^{\#}\right) $ explicitly in the standard PBW
basis and this basis is dual to the exponential coordinates for the group $%
\mathfrak{G}^{\#}$.

To derive the transformation $\mathfrak{g}\longrightarrow \mathfrak{g}^{\#}$
one must identify the elements in two Hopf algebras $U\left( \mathfrak{g}%
^{\#}\right) $ and $U\left( \mathfrak{g}_{\left(
\mathfrak{g}\right) }^{\#}\right) $%
. The relations in them are induced by the ideals $I$ that encode the Lie
multiplications. It is sufficient to compare the ideals $I\left( \mathfrak{g}%
^{\#}\right) $ and
$I\left( \mathfrak{g}_{\left( \mathfrak{g}\right) }^{\#}\right) $%
. This gives the finite set of relations
\begin{equation}
e_{i}^{\#\ast }=\phi _{i}^{\ast }\left( \left\{ e_{j}^{\ast }\right\} \right)
\label{dual-rel}
\end{equation}
that describe $\mathfrak{g}^{\#}$-basic
elements in terms of $\mathfrak{g}$-basis
and thus introduces the $\mathfrak{g}^{\#}$-PBW coordinates
in $\left( \mathrm{Fun%
}\left( \mathfrak{G}^{\#}\right) \right) ^{\ast }$. Constructing the relations
dual to (\ref{dual-rel}) we return to $\mathrm{Fun}\left( \mathfrak{G}%
^{\#}\right) $ and find the desired transformation
\begin{equation*}
e_{i}^{\#}=\phi _{i}\left( \left\{ e_{j}\right\} \right) .
\end{equation*}

Finally this transformation is to be applied to the deformed algebra $U_{%
\mathcal{F}}\left( \mathfrak{g}\right) =$ $\mathrm{Fun}_{\mathrm{def}}\left(
\mathfrak{G}^{\#}\right) $. This gives the deformed group $\mathfrak{G}^{\#}$
in its exponential coordinates. In this presentation the comultiplication of
functions on the solvable group $\mathfrak{G}^{\#}$ become highly transparent
demonstrating the exponential map of the adjoint action of the dual Lie
algebra $\mathfrak{g}^{\#}$.

\section{Example. Twisted algebra $U_{\mathcal{F}}\left( sl\left( 3\right)
\right) $ in dual group coordinates}

Let $\mathfrak{g}=\mathfrak{sl}\left( 3\right) $ and
consider the $r$-elements: $%
r_{J}\left( \gamma \right) =h \wedge e_{13}$ and $%
r_{E}=e_{12}\wedge e_{23},$ where $h =h_{13}+\gamma
h_{\perp }$. Here the generators $e_{ab}$ are the ordinary matrix units
and the Cartan elements are $h_{ab}=1/2 \left( e_{aa} - e_{bb} \right)$,
$h_{\perp }=1/2 \left( e_{11} - 2e_{22}+ e_{33} \right)$.
The sum $r_{EJ}\left( \gamma \right) =r_{J}\left( \gamma
\right) +r_{E}$ is the set of solutions of the CYBE. Or equivalently,
the sum of deforming functions $f_{J}^{\#}\left(
\gamma \right) +f_{E}^{\#}$ defines the Lie multiplication $\left[ ,\right]
^{\#}:=\left( f_{J}^{\#}\left( \gamma \right) +f_{E}^{\#}\right) \left(
,\right) $. This bracket is the Lie composition of the
dual algebra $\mathfrak{sl}(3)^{\#}$ on the
space $V_{\mathfrak{g}}^{\ast }$. The basic commutation relations are
\begin{equation}
\begin{array}{ll}
\left[ e_{13}^{\ast },h^{\ast }\right] ^{\#}=h^{\ast }, & \left[ h^{\ast
},e_{12}^{\ast }\right] ^{\#}=\left[ h^{\ast
},e_{23}^{\ast }\right] ^{\#}=0, \\
\left[ e_{13}^{\ast },e_{21}^{\ast }\right] ^{\#}=+\frac{1}{2}\left( 3\gamma
+1\right) e_{21}^{\ast }, & \left[ h^{\ast },h_{\perp
}^{\ast }\right] ^{\#}=-2\gamma e_{31}^{\ast }, \\
\left[ e_{13}^{\ast },e_{12}^{\ast }\right] ^{\#}=-\frac{1}{2}\left( 3\gamma
-1\right) e_{12}^{\ast }, & \left[ e_{13}^{\ast },e_{32}^{\ast }\right]
^{\#}=-\frac{1}{2}\left( 3\gamma -1\right) e_{32}^{\ast }, \\
\begin{array}{l}
\left[ e_{13}^{\ast },e_{31}^{\ast }\right] ^{\#}=\left[ e_{23}^{\ast
},e_{32}^{\ast }\right] ^{\#}= \\
=\left[ e_{12}^{\ast },e_{21}^{\ast }\right] ^{\#}=e_{31}^{\ast },
\end{array}
& \left[ e_{13}^{\ast },e_{23}^{\ast }\right] ^{\#}=+\frac{1}{2}\left(
3\gamma +1\right) e_{23}^{\ast }, \\
\left[ h_{\perp }^{\ast },e_{23}^{\ast }\right] ^{\#}=\left( 1-\gamma
\right) e_{21}^{\ast }, & \left[ e_{12}^{\ast },e_{23}^{\ast }\right]
^{\#}=-h\left( \gamma \right) ^{\ast }. \\
\left[ h_{\perp }^{\ast },e_{12}^{\ast }\right] ^{\#}=\left( 1+\gamma
\right) e_{32}^{\ast }, &
\end{array}
\label{dual-sl3-alg}
\end{equation}
This algebra has the four-dimensional subalgebra (equivalent to the
subalgebra $\mathfrak{g}_{c}$, the carrier of the $r$-matrix $r_{EJ}\left(
\gamma \right) $) with the weight system $\Gamma _{c}^{\#}=\left\{ 0,\varphi
_{31},\varphi _{32},\varphi _{21}\right\} $ and the four-dimensional Abelian
ideal generated by the set $\left\{ h_{\perp }^{\ast },e_{31}^{\ast
},e_{32}^{\ast },e_{21}^{\ast }\right\} $ with the weights $\Gamma
_{a}^{\#}=\left\{ 0,\rho _{31},\rho _{32},\rho _{21}\right\} $ . Here $%
\left\{ \varphi _{ab}\right\} $ and $\left\{ \rho _{ab}\right\} $ are the
copies of those vectors from $\Lambda \left( \mathfrak{sl}\left( 3\right)
\right) $ that correspond to the generators $e_{ab}$ . In this particular
case both subsystems contain the same sets of vectors. The action of $%
\mathrm{ad}\left( \mathfrak{g}_{c}^{\#}\right) $ on $h_{\perp }^{\ast }$ shows
that $h_{\perp }^{\ast }$ behaves as ''orthogonal'' to $h^{\ast }$.

The following twisting element
\begin{equation}
\mathcal{F}\left( \gamma \right) =\mathcal{F}_{E}\left( \gamma \right)
\mathcal{F}_{J}\left( \gamma \right) =\exp \left( \xi e_{12}\otimes e_{23}e^{%
\frac{1}{2}\left( 3\gamma -1\right) \sigma \left( \xi \right) }\right) \exp
\left( h\otimes \sigma \left( \xi \right) \right)  \label{sl3-tw}
\end{equation}
is the solution of the twist equations \cite{Drin-1} corresponding to the $r$%
-matrix $r_{EJ}\left( \gamma \right) $. Here $\sigma \left( \xi \right) =\ln
\left( 1+\xi e_{13}\right) $.

The parameter $\gamma $ describes the Reshetikhin rotation of the Jordanian
factor while $\xi $ is the deformation parameter. The latter corresponds to
the automorphism
\begin{equation*}
e_{\tau }\longrightarrow \xi ^{\tau \left( h_{13}\right) }\widetilde{e_{\tau
}},
\end{equation*}
where $\tau \in \Gamma \left( \mathfrak{sl}\left( 3\right) \right) $ are
the weights and the element $h_{13}$ performs the gradation in $%
\Gamma $, $h_{13}:\Gamma \longrightarrow $ $R^{1}$. Together with the
scaling this gives the parametrization appropriate to perform the second
classical limit
\begin{equation}
e_{\tau }\longrightarrow \frac{\xi ^{\tau \left( h_{13}\right) }}{%
\varepsilon }\widehat{e_{\tau }},\qquad \frac{\xi }{\varepsilon }=\zeta .
\label{sl3-scale}
\end{equation}

The twisting is the transformation of the comultiplication
in $U\left( \mathfrak{g}\right) $ performed by the operator
\begin{equation*}
\prod_{q=p}^{1}e^{\mathrm{ad}
\left( \ln \mathcal{F}_{q}\right) }=e^{\frac{1}{%
\varepsilon }\mathrm{ad}
\left( \mathrm{BCH}\left\{ \Psi _{\left( \mathcal{F}%
\right) q}\left( \left\{ \widehat{e_{i}};\zeta _{s}\right\} \right)
;\varepsilon \right\} \right) }
\end{equation*}
Here $\mathrm{BCH}$ denotes the Baker-Campbell-Hausdorff series and $\Psi
_{\left( \mathcal{F}\right) q}\left( \left\{ \widehat{e_{i}};\zeta
_{s}\right\} \right) $ are the logarithms of the twisting factors $\mathcal{F%
}_{j}$ :
\begin{equation}
\Psi _{\left( \mathcal{F}\right) q}\left( \left\{ \widehat{e_{i}};\zeta
_{s}\right\} \right) =\varepsilon \ln \mathcal{F}_{q},\qquad \Psi \in
U\otimes U;\quad q=1,\ldots ,p;i=1,\ldots ,n.
\end{equation}
For the twisting element (\ref{sl3-tw}) these are
\begin{eqnarray*}
\Psi _{\left( \mathcal{F}\right) 1}\left( \left\{ \widehat{e_{i}};\zeta
\right\} \right) &=&\widehat{h}\otimes \widehat{\sigma }\left( \zeta \right)
, \\
\Psi _{\left( \mathcal{F}\right) 2}\left( \left\{ \widehat{e_{i}};\zeta
\right\} \right) &=&\zeta \widehat{e_{12}}\otimes \widehat{e_{23}}e^{\frac{1%
}{2}\left( 3\gamma -1\right) \widehat{\sigma }\left( \zeta \right) }.
\end{eqnarray*}
In the second classical limit only the zero power terms in the $\mathrm{BCH}$%
-series remain and the compositions of
the dual group $\mathfrak{G}^{\#}=\left(
\mathfrak{SL}(3)\right) ^{\#}$ in $\mathfrak{g}=\mathfrak{sl}(3)$-coordinates,
i.e. the coproducts
in $\mathrm{Fun}\left( \mathfrak{G}^{\#}\right) $ in $\mathfrak{g}$%
-basis, can be obtained by the formula:
\begin{eqnarray*}
\Delta _{\mathcal{F},\left\{ \zeta \right\} }^{\lim }\left( \widehat{e_{j}}%
\right) &=&\lim_{\varepsilon \longrightarrow 0}\left( e^{\frac{1}{%
\varepsilon }\mathrm{ad}\left( \mathrm{BCH}\left\{ \Psi _{\left( \mathcal{F}%
\right) q}\left( \left\{ \widehat{e_{i}};\zeta _{s}\right\} \right)
;\varepsilon \right\} \right) }\circ \left( \widehat{e_{j}}\otimes
1+1\otimes \widehat{e_{j}}\right) \right) = \\
&=&\lim_{\varepsilon \longrightarrow 0}\left( e^{\frac{1}{\varepsilon }%
\mathrm{ad}\left( \Psi _{\left( \mathcal{F}\right) 1}\left( \left\{ \widehat{%
e_{i}};\zeta \right\} \right) +\Psi _{\left( \mathcal{F}\right) 2}\left(
\left\{ \widehat{e_{i}};\zeta \right\} \right) \right) }\circ \left(
\widehat{e_{j}}\otimes 1+1\otimes \widehat{e_{j}}\right) \right)
\end{eqnarray*}
In particular the scaled $\mathfrak{g}$-basic elements
(and $\widehat{\sigma }%
\left( \zeta \right) $) have the following coproducts in $\mathrm{Fun}\left(
\mathfrak{G}^{\#}\right) $ :
\begin{equation}
\begin{array}{lll}
\Delta _{\mathcal{F}}^{\lim }(\widehat{h}\left( \gamma \right) ) & = &
\widehat{h}\left( \gamma \right) \otimes e^{-\widehat{\sigma }\left( \zeta
\right) }+1\otimes \widehat{h}\left( \gamma \right) -\xi \widehat{e_{12}}%
\otimes \widehat{e_{23}}e^{\frac{3}{2}\left( \gamma -1\right) \widehat{%
\sigma }\left( \zeta \right) }; \\
\Delta _{\mathcal{F}}^{\lim }(\widehat{h_{\perp }}) & = & \widehat{h_{\perp }%
}\otimes 1+1\otimes \widehat{h_{\perp }}; \\
\Delta _{\mathcal{F}}^{\lim }(\widehat{e_{12}}) & = & \widehat{e_{12}}%
\otimes e^{^{\frac{1}{2}\left( 3\gamma -1\right) \widehat{\sigma }\left(
\zeta \right) }}+1\otimes \widehat{e_{12}}; \\
\Delta _{\mathcal{F}}^{\lim }(\widehat{e_{23}}) & = & \widehat{e_{23}}%
\otimes e^{-\frac{1}{2}\left( 3\gamma -1\right) \widehat{\sigma }\left(
\zeta \right) }+e^{\widehat{\sigma }\left( \zeta \right) }\otimes e\widehat{%
_{23}}; \\
\Delta _{\mathcal{F}}^{\lim }(\widehat{\sigma }\left( \zeta \right) ) & = &
\widehat{\sigma }\left( \zeta \right) \otimes 1+1\otimes \widehat{\sigma }%
\left( \zeta \right) ; \\
\Delta _{\mathcal{F}}^{\lim }(\widehat{e_{21}}) & = & \widehat{e_{21}}%
\otimes e^{-\frac{1}{2}\left( 3\gamma +1\right) \widehat{\sigma }\left(
\zeta \right) }+1\otimes \widehat{e_{21}}+\zeta \left( 1-\gamma \right)
\widehat{h_{\perp }}\otimes \widehat{e_{23}}e^{-\widehat{\sigma }\left(
\zeta \right) }; \\
\Delta _{\mathcal{F}}^{\lim }(\widehat{e_{32}}) & = & \widehat{e_{32}}%
\otimes e^{^{\frac{1}{2}\left( 3\gamma -1\right) \widehat{\sigma }\left(
\zeta \right) }}+1\otimes \widehat{e_{32}}+\zeta \widehat{h}\left( \gamma
\right) \otimes \widehat{e_{12}}e^{-\widehat{\sigma }\left( \zeta \right) }+
\\
&  &
\begin{array}{l}
+\zeta \widehat{e_{12}}\otimes \left( \widehat{h}\left( \gamma \right)
-\left( \gamma +1\right) \widehat{h_{\perp }}\right) e^{^{\frac{1}{2}\left(
3\gamma -1\right) \widehat{\sigma }\left( \zeta \right) }}- \\
-\zeta \widehat{h}\left( \gamma \right) \widehat{e_{12}}\otimes \left( e^{^{%
\frac{1}{2}\left( 3\gamma -1\right) \widehat{\sigma }\left( \zeta \right)
}}-e^{^{\frac{3}{2}\left( \gamma -1\right) \widehat{\sigma }\left( \zeta
\right) }}\right) - \\
-\zeta ^{2}\widehat{e_{12}}\otimes \widehat{e_{23}}\widehat{e_{12}}e^{^{%
\frac{3}{2}\left( \gamma -1\right) \widehat{\sigma }\left( \zeta \right)
}}-\zeta ^{2}\widehat{e_{12}}^{2}\otimes \widehat{e_{23}}e^{\left( 3\gamma
-2\right) \widehat{\sigma }\left( \zeta \right) };
\end{array}
\\
\Delta _{\mathcal{F}}^{\lim }(\widehat{e_{31}}) & = & \widehat{e_{31}}%
\otimes e^{-\widehat{\sigma }\left( \zeta \right) }+1\otimes \widehat{e_{31}}%
+ \\
&  &
\begin{array}{l}
+2\zeta \widehat{h}\left( \gamma \right) \otimes \left( \widehat{h}\left(
\gamma \right) -\gamma \widehat{h_{\perp }}\right) e^{-\widehat{\sigma }%
\left( \zeta \right) } \\
-\zeta \widehat{h}\left( \gamma \right) \left( \widehat{h}\left( \gamma
\right) -2\gamma \widehat{h_{\perp }}\right) \otimes \left( e^{-\widehat{%
\sigma }\left( \zeta \right) }-e^{-2\widehat{\sigma }\left( \zeta \right)
}\right) - \\
+\zeta ^{2}\widehat{h}\left( \gamma \right) \widehat{e_{12}}\otimes \widehat{%
e_{23}}\left( e^{^{\frac{3}{2}\left( \gamma -1\right) \widehat{\sigma }%
\left( \zeta \right) }}-2e^{^{\frac{1}{2}\left( 3\gamma -5\right) \widehat{%
\sigma }\left( \zeta \right) }}\right) + \\
+\zeta \widehat{e_{12}}\otimes \widehat{e_{21}}e^{^{\frac{1}{2}\left(
3\gamma -1\right) \widehat{\sigma }\left( \zeta \right) }}-\zeta \widehat{%
e_{32}}\otimes \widehat{e_{23}}e^{^{\frac{3}{2}\left( \gamma -1\right)
\widehat{\sigma }\left( \zeta \right) }}- \\
-2\zeta ^{2}\widehat{e_{12}}\otimes \left( \widehat{h}\left( \gamma \right)
-\gamma \widehat{h_{\perp }}\right) \widehat{e_{23}}e^{^{\frac{3}{2}\left(
\gamma -1\right) \widehat{\sigma }\left( \zeta \right) }} \\
+\zeta \widehat{e_{12}}^{2}\otimes \widehat{e_{23}}^{2}e^{^{3\left( \gamma
-1\right) \widehat{\sigma }\left( \zeta \right) }};
\end{array}
\end{array}
\label{lim-sl3-group}
\end{equation}

The standard dualization transforms this costructure into the multiplication
for the Hopf algebra $\left( U_{\mathcal{F}}\left( \mathfrak{sl}\left( 3\right)
\right) \right) ^{\ast }=\left( \mathrm{Fun}\left( \mathfrak{SL}\left( 3\right)
^{\#}\right) \right) ^{\ast }
=U\left( \left( \mathfrak{sl}(3)\right) _{\mathfrak{g}%
}^{\#}\right) $. The latter is thus obtained in terms of $\mathfrak{g}$%
-coordinates $\left\{ \widehat{e}_{\tau }^{\ast }\right\} $. The same Hopf
algebra in terms of $\mathfrak{g}^{\#}$-coordinates $\left\{ e_{\tau }^{\#\ast
}\right\} $ is defined by the relations (\ref{dual-sl3-alg}) that generate
the necessary ideal in the tensor algebra over $V_{\mathfrak{g}^{\#}}$. The
deformation parameter $\zeta $ is introduced by the scaling (\ref{sl3-scale}%
). Comparing the associative multiplications in $U\left( \left( \mathfrak{sl}%
(3)\right) ^{\#}\right) $ and
$U\left( \left( \mathfrak{sl}(3)\right) _{\mathfrak{g}%
}^{\#}\right) $ we find two types of nontrivial relations. The relations of
the first type are
\begin{equation*}
\begin{array}{ll}
\widehat{\sigma \left( \zeta \right) }^{\ast }\cdot \widehat{h\left( \gamma
\right) }^{\ast }=\widehat{h\left( \gamma \right) }^{\ast };\ldots &
e_{13}^{\#\ast }\left( \zeta \right) \cdot h\left( \gamma \right) ^{\#\ast
}=h\left( \gamma \right) ^{\#\ast };\ldots
\end{array}
\end{equation*}
and
\begin{equation*}
\begin{array}{ll}
\widehat{\sigma \left( \zeta \right) }^{\ast }\cdot \left( \widehat{e_{23}}%
e^{-\widehat{\sigma \left( \zeta \right) }}\right) ^{\ast }=\frac{1}{2}%
\left( \widehat{e_{23}}e^{-\widehat{\sigma \left( \zeta \right) }}\right)
^{\ast }; & e_{13}^{\#\ast }\left( \zeta \right) \cdot e_{23}^{\#\ast }=%
\frac{1}{2}e_{23}^{\#\ast }.
\end{array}
\end{equation*}
They signify that the transition to $\mathfrak{g}^{\#}$- generators
includes the substitutions
\begin{eqnarray*}
e_{13}^{\#}\left( \zeta \right) &=&\ln \left( 1+\zeta \widehat{e_{13}}%
\right) , \\
e_{23}^{\#} &=&\widehat{e_{23}}e^{-\widehat{\sigma \left( \zeta \right) }}.
\end{eqnarray*}
These transformations are defined in the space $U\left( \mathfrak{g}%
_{c}^{\#}\right) $ of the carrier subalgebra. They coincide with those induced
by the twisting element morphism $\mathcal{F}:%
\mathfrak{g}_{c}\longrightarrow \mathfrak{g}_{c}^{\#}$.

Much more important are the relations on the space $V_{\mathfrak{a}^{\#}}$
complimentary to $V_{\mathfrak{g}_{c}^{\#}}$. They cannot be extracted from the
twisting element itself and depend on the subrepresentation $\mathrm{ad}%
\mathfrak{g}_{c}^{\#}|_{V_{\mathfrak{a}^{\#}}}$. In our case there are only two
relations of this type obtained from (\ref{lim-sl3-group}) and (\ref
{dual-sl3-alg}):
\begin{equation*}
\begin{array}{ll}
\widehat{h\left( \gamma \right) }^{\ast }\cdot \widehat{e}_{12}^{\ast
}=\left( \widehat{h\left( \gamma \right) }\cdot \widehat{e}_{12}\right)
^{\ast }+\zeta \widehat{e}_{12}^{\ast }; & h\left( \gamma \right) ^{\#\ast
}\cdot e_{12}^{\#\ast }=\left( h\left( \gamma \right) ^{\#}\cdot
e_{12}^{\#}\right) ^{\ast }; \\
\begin{array}{l}
\widehat{h\left( \gamma \right) }^{\ast }\cdot \widehat{h\left( \gamma
\right) }^{\ast }= \\
=2\left( \left( \widehat{h\left( \gamma \right) }\right) ^{2}\right) ^{\ast
}+2\zeta \widehat{e}_{31}^{\ast };
\end{array}
&
\begin{array}{l}
h\left( \gamma \right) ^{\#\ast }\cdot h\left( \gamma \right) ^{\#\ast }= \\
=2\left( \left( h\left( \gamma \right) ^{\#}\right) ^{2}\right) ^{\ast };
\end{array}
\end{array}
\end{equation*}
This gives two nontrivial relations
between $\left( \mathfrak{g}\right) ^{\ast }$%
- and $\left( \mathfrak{g}^{\#}\right) ^{\ast }$- generators
\begin{eqnarray*}
\left( \widehat{h\left( \gamma \right) }\cdot \widehat{e}_{12}\right) ^{\ast
}\left( \zeta \right) &=&\left( h\left( \gamma \right) ^{\#}\cdot
e_{12}^{\#}\right) ^{\ast }-\zeta e_{32}^{\#\ast }, \\
\left( \left( \widehat{h\left( \gamma \right) }\right) ^{2}\right) ^{\ast
}\left( \zeta \right) &=&\frac{1}{2}\left( \left( h\left( \gamma \right)
^{\#}\right) ^{2}\right) ^{\ast }-\zeta e_{31}^{\#\ast },
\end{eqnarray*}
and correspondingly between $\mathfrak{g}$- and
 $\mathfrak{g}^{\#}$- generators
\begin{eqnarray*}
e_{32}^{\#}\left( \zeta \right) &=&\widehat{e}_{32}-\zeta \widehat{h\left(
\gamma \right) }\widehat{e}_{12}, \\
e_{31}^{\#}\left( \zeta \right) &=&\widehat{e}_{31}-\zeta \left( \widehat{%
h\left( \gamma \right) }\right) ^{2}.
\end{eqnarray*}
As we have already mentioned the basic
elements $\left\{ e_{\tau }^{\#}\right\} $
constructed above in terms of the undeformed group coordinates $\widehat{e}%
_{\tau }$ are appropriate also for the twisted Hopf algebra. Summing up we
get the $\mathfrak{g}\longrightarrow \mathfrak{g}^{\#}$ transformation:
\begin{equation}
\begin{array}{rcl}
e_{13}^{\#}\left( \zeta \right) &=&\sigma \left( \zeta \right) ,\\
e_{23}^{\#}\left( \zeta \right) &=&e_{23}e^{-\sigma \left( \zeta \right) },\\
e_{32}^{\#}\left( \zeta \right) &=&e_{32}-\zeta h\left( \gamma \right)
e_{12},  \notag \\
e_{31}^{\#}\left( \zeta \right) &=&e_{31}-\zeta h\left( \gamma \right) ^{2}.
\end{array}
\label{dual-sl3-coord}
\end{equation}

In terms of the original $\mathfrak{g}$-coordinates (and the function $\sigma
\left( \xi \right) =\ln \left( 1+\xi e_{13}\right) $) the costructure in the
twisted
algebra $U_{\mathcal{F}}\left( \mathfrak{sl}\left( 3\right) \right) $ is
defined by the following relations \cite{KLM},\cite{LO}:
\begin{equation}
\begin{array}{lll}
\Delta _{\mathcal{F}}(h\left( \gamma \right) ) & = & h\left( \gamma \right)
\otimes e^{-\sigma \left( \xi \right) }+1\otimes h\left( \gamma \right) -\xi
e_{12}\otimes e_{23}e^{\frac{3}{2}\left( \gamma -1\right) \sigma \left( \xi
\right) }; \\
\Delta _{\mathcal{F}}(h_{\perp }) & = & h_{\perp }\otimes 1+1\otimes
h_{\perp }; \\
\Delta _{\mathcal{F}}(e_{12}) & = & e_{12}\otimes e^{^{\frac{1}{2}\left(
3\gamma -1\right) \sigma \left( \xi \right) }}+1\otimes e_{12}; \\
\Delta _{\mathcal{F}}(e_{23}) & = & e_{23}\otimes e^{-\frac{1}{2}\left(
3\gamma -1\right) \sigma \left( \xi \right) }+e^{\sigma \left( \xi \right)
}\otimes e_{23}; \\
\Delta _{\mathcal{F}}(\sigma \left( \xi \right) ) & = & \sigma \left( \xi
\right) \otimes 1+1\otimes \sigma \left( \xi \right) ; \\
\Delta _{\mathcal{F}}(e_{21}) & = & e_{21}\otimes e^{-\frac{1}{2}\left(
3\gamma +1\right) \sigma \left( \xi \right) }+1\otimes e_{21}+\xi \left(
1-\gamma \right) h_{\perp }\otimes e_{23}e^{-\sigma \left( \xi \right) }; \\
\Delta _{\mathcal{F}}(e_{32}) & = & e_{32}\otimes e^{^{\frac{1}{2}\left(
3\gamma -1\right) \sigma \left( \xi \right) }}+1\otimes e_{32}+\xi h\left(
\gamma \right) \otimes e_{12}e^{-\sigma \left( \xi \right) }+ \\
&  &
\begin{array}{l}
+\xi e_{12}\otimes \left( h\left( \gamma \right) -\left( \gamma +1\right)
h_{\perp }\right) e^{^{\frac{1}{2}\left( 3\gamma -1\right) \sigma \left( \xi
\right) }}- \\
-\xi h\left( \gamma \right) e_{12}\otimes \left( e^{^{\frac{1}{2}\left(
3\gamma -1\right) \sigma \left( \xi \right) }}-e^{^{\frac{3}{2}\left( \gamma
-1\right) \sigma \left( \xi \right) }}\right) - \\
-\xi ^{2}e_{12}\otimes e_{23}e_{12}e^{^{\frac{3}{2}\left( \gamma -1\right)
\sigma \left( \xi \right) }}-\xi ^{2}e_{12}^{2}\otimes e_{23}e^{\left(
3\gamma -2\right) \sigma \left( \xi \right) };
\end{array}
\\
\Delta _{\mathcal{F}}(e_{31}) & = & e_{31}\otimes e^{-\sigma \left( \xi
\right) }+1\otimes e_{31}+ \\
&  &
\begin{array}{l}
+2\xi h\left( \gamma \right) \otimes \left( h\left( \gamma \right) -\gamma
h_{\perp }\right) {e^{-\sigma \left( \xi \right) }} \\
+\xi \left( h\left( \gamma \right) -h\left( \gamma \right) ^{2}\right)
\otimes \left( {e^{-\sigma \left( \xi \right) }-e^{-2\,\sigma \left( \xi
\right) }}\right) \\
+\xi e_{{1,2}}\otimes e_{{2,1}}{e^{\frac{1}{2}\,\left( 3\,\gamma -1\right)
\sigma \left( \xi \right) }} \\
-\xi e_{{3,2}}\otimes e_{{2,3}}{e^{\frac{3}{2}\,\left( \gamma -1\right)
\sigma \left( \xi \right) }}+\xi ^{2}h\left( \gamma \right) e_{{1,2}}\otimes
e_{{2,3}}{e^{\frac{3}{2}\,\left( \gamma -1\right) \sigma \left( \xi \right) }%
} \\
+2\xi ^{2}e_{{1,2}}\otimes e_{{2,3}}{e^{\frac{1}{2}\,\left( 3\,\gamma
-5\right) \sigma \left( \xi \right) }}-\xi ^{2}e_{{1,2}}\otimes e_{{2,3}}{e^{%
\frac{3}{2}\,\left( \gamma -1\right) \sigma \left( \xi \right) }} \\
-2\xi ^{2}e_{{1,2}}\otimes \left( h\left( \gamma \right) -\gamma h_{\perp
}\right) e_{{2,3}}{e^{\frac{3}{2}\,\left( \gamma -1\right) \sigma \left( \xi
\right) }} \\
-2\xi ^{2}h\left( \gamma \right) e_{{1,2}}\otimes e_{{2,3}}{e^{\frac{1}{2}%
\,\left( 3\,\gamma -5\right) \sigma \left( \xi \right) }}+\xi ^{3}e_{{1,2}%
}^{2}\otimes e_{{2,3}}^{2}e^{3\,\left( \gamma -1\right) \sigma \left( \xi
\right) }
\end{array}
\end{array}
\end{equation}
Applying to these coproducts the transformation (\ref{dual-sl3-coord}) we
get the same costructure in $\mathfrak{g}^{\#}$-coordinates:
\begin{equation}
\begin{array}{lll}
\Delta _{\mathcal{F}}(h^{\#}\left( \gamma \right) ) & = & h^{\#}\left(
\gamma \right) \otimes e^{-e_{13}^{\#}}+1\otimes h^{\#}\left( \gamma \right)
-\zeta e_{12}^{\#}\otimes e_{23}^{\#}e^{\frac{1}{2}\left( \gamma -1\right)
e_{13}^{\#}}; \\
\Delta _{\mathcal{F}}(e_{13}^{\#}) & = & e_{13}^{\#}\otimes 1+1\otimes
e_{13}^{\#}; \\
\Delta _{\mathcal{F}}(e_{12}^{\#}) & = & e_{12}^{\#}\otimes e^{^{\frac{1}{2}%
\left( 3\gamma -1\right) e_{13}^{\#}}}+1\otimes e_{12}^{\#}; \\
\Delta _{\mathcal{F}}(e_{23}^{\#}) & = & e_{23}^{\#}\otimes e^{-\frac{1}{2}%
\left( 3\gamma +1\right) e_{13}^{\#}}+1\otimes e_{23}^{\#}; \\
\Delta _{\mathcal{F}}(h_{\perp }^{\#}) & = & h_{\perp }^{\#}\otimes
1+1\otimes h_{\perp }^{\#}; \\
\Delta _{\mathcal{F}}(e_{21}^{\#}) & = & e_{21}^{\#}\otimes e^{-\frac{1}{2}%
\left( 3\gamma +1\right) e_{13}^{\#}}+1\otimes e_{21}^{\#}+\zeta \left(
1-\gamma \right) h_{\perp }^{\#}\otimes e_{23}^{\#}; \\
\Delta _{\mathcal{F}}(e_{32}^{\#}) & = & e_{32}^{\#}\otimes e^{^{\frac{1}{2}%
\left( 3\gamma -1\right) e_{13}^{\#}}}+1\otimes e_{32}^{\#}-\zeta \left(
\gamma +1\right) e_{12}^{\#}\otimes h_{\perp }^{\#}e^{^{\frac{1}{2}\left(
3\gamma -1\right) e_{13}^{\#}}}; \\
\Delta _{\mathcal{F}}(e_{31}^{\#}) & = & e_{31}^{\#}\otimes
e^{-e_{13}^{\#}}+1\otimes e_{31}^{\#}+ \\
&  &
\begin{array}{l}
+\zeta e_{12}^{\#}\otimes e_{21}^{\#}e{^{\frac{1}{2}\,\left( 3\,\gamma
-1\right) e_{13}^{\#}}}-\zeta e_{32}^{\#}\otimes e_{23}^{\#}e{^{\frac{3}{2}%
\,\left( \gamma -1\right) e_{13}^{\#}}} \\
-2\zeta \gamma h^{\#}\left( \gamma \right) \otimes h_{\perp }^{\#}e{%
^{-e_{13}^{\#}}}+2\zeta ^{2}\gamma e_{12}^{\#}\otimes h_{\perp
}^{\#}e_{23}^{\#}e{^{\frac{3}{2}\,\left( \gamma -1\right) e_{13}^{\#}};}
\end{array}
\end{array}
\label{costrsl3-in-duals}
\end{equation}
This presentation reveals the $\mathfrak{G}^{\#}$-group law in $U_{\mathcal{F}%
}\left( \mathfrak{sl}\left( 3\right) \right) $ described in terms of the
exponential coordinates. The first four relations correspond to the group
multiplication in $\mathfrak{G}_{c}^{\#}$, the other four expose the adjoint
action of $\mathfrak{G}_{c}^{\#}$ on the 4-dimensional
space $V_{\mathfrak{a}%
^{\#}}$. In the standard orthonormal basis $\left\{ \mathbf{e}_{i}\right\} $
the weights of the diagram
\begin{equation*}
\Gamma _{\mathcal{F}}^{\#}=\Gamma _{c}^{\#}\cup \Gamma _{a}^{\#}=\left\{
0,\varphi _{31},\varphi _{32},\varphi _{21}\right\} \cup \left\{ 0,\rho
_{31},\rho _{32},\rho _{21}\right\}
\end{equation*}
have the form:
\begin{equation*}
\varphi _{ab}=\rho _{ab}=\mathbf{e}_{a}-\mathbf{e}_{b}.
\end{equation*}
Consider, for example, the coproduct $\Delta _{\mathcal{F}}(e_{31}^{\#})$.
The first two expressions are due to the trivial property of the unity and
the adjoint action
\begin{equation*}
\left[ e_{13}^{\#\ast },e_{31}^{\#\ast }\right] =\ldots +e_{31}^{\#\ast
}+\ldots
\end{equation*}
The remaining four sets of terms are directly correlated with the weights
shifts:
\begin{equation*}
\left.
\begin{array}{c}
\varphi _{32}\circ \rho _{21} \\
\varphi _{21}\circ \rho _{32} \\
\varphi _{31}\circ \rho \left( h_{\perp }^{\#}\right) \\
\varphi _{32}\circ \varphi _{21}\circ \rho \left( h_{\perp }^{\#}\right)
\end{array}
\right\} =\rho _{31}.
\end{equation*}

The redeveloped costructure (\ref{costrsl3-in-duals}) expose explicitly the
properties that can be used to find new twist cocycles. First of all it is
clearly seen that the parametrized set of dual groups $\mathfrak{G}^{\#}\left(
\gamma \right) $ has tree irregular points,
\begin{equation*}
\gamma =0,\pm 1.
\end{equation*}
At the point $\gamma =0$ the coproducts for $e_{31}^{\#}$ has no terms
containing $h_{\perp }^{\#}$. The shifts $\varphi _{31}\circ \rho \left(
h_{\perp }^{\#}\right) ,$ and $\varphi _{32}\circ \varphi _{21}\circ \rho
\left( h_{\perp }^{\#}\right) $ result in zeros because in this case $h^{\#}$
is orthogonal to the weight of $h_{\perp }^{\#}$ .

The points $\gamma =+1$ $\left( \gamma =-1\right) $\ are especially
interesting. In them the coproduct $\Delta _{\mathcal{F}}(e_{32}^{\#})$
(correspondingly $\Delta _{\mathcal{F}}(e_{21}^{\#})$) becomes
quasiprimitive. Together with the primitivity of $\Delta _{\mathcal{F}%
}(h_{\perp }^{\#})$ this means that in these points the twist equations for
the Hopf algebra $U_{\mathcal{F}}\left( \mathfrak{sl}\left(
3\right) \right) $ have the additional solutions:
\begin{eqnarray*}
\mathcal{F}_{\mathcal{J}_{P+}} &=&\exp \left( -\frac{2}{3}h_{\perp
}^{\#}\otimes \ln \left( \left( 1+e_{21}^{\#}\right) e^{2e_{13}^{\#}}\right)
\right) \qquad \mathrm{for}\quad \gamma =+1, \\
\mathcal{F}_{\mathcal{J}_{P-}} &=&\exp \left( +\frac{2}{3}h_{\perp
}^{\#}\otimes \ln \left( \left( 1+e_{32}^{\#}\right) e^{2e_{13}^{\#}}\right)
\right) \\
&=&\exp \left( +\frac{2}{3}h_{\perp }\otimes \ln \left( \left(
1+e_{32}-\zeta he_{12}\right) e^{2\sigma _{13}\left( \zeta \right) }\right)
\right) \qquad \mathrm{for}\quad \gamma =-1.
\end{eqnarray*}
Correspondingly for $U\left( \mathfrak{sl}\left( 3\right) \right) $ we have two
parabolic twists :
\begin{eqnarray*}
\mathcal{F}_{\mathfrak{P}+} &=&\mathcal{F}_{\mathcal{J}_{P+}}\mathcal{F}_{%
\mathcal{E}}\left( \gamma \right) \mathcal{F}_{\mathcal{J}}\left( \gamma \right)
\qquad \mathrm{for}\quad \gamma =+1, \\
\mathcal{F}_{\mathfrak{P}-} &=&\mathcal{F}_{\mathcal{J}_{P-}}\mathcal{F}_{%
\mathcal{E}}\left( \gamma \right) \mathcal{F}_{\mathcal{J}}\left( \gamma \right)
\qquad \mathrm{for}\quad \gamma =-1.
\end{eqnarray*}
This construction generalizes the elementary parabolic twist first presented
in \cite{Lya-Sam}.

Both twists $\mathcal{F}_{\mathcal{J}_{P+}}$ and $\mathcal{F}_{\mathcal{J}%
_{P-}}$ are the deformed Jordanian factors and the canonical form $\mathcal{F%
}_{\mathcal{J}}=e^{H\otimes \sigma }$ is reobtained when $\zeta =0$. This
means that when the preceding extended Jordanian twist is trivialized, $%
\mathcal{F}_{\mathcal{E}}\mathcal{F}_{\mathcal{J}}|_{\zeta =0}=1\otimes 1$ ,
the jordanian deformation $\mathcal{F}_{\mathcal{J}}=e^{\mp \frac{2}{3}%
h_{\perp }\otimes \sigma _{\binom{21}{32}}}$ remains possible with the
ordinary function $\sigma _{\binom{21}{32}}=\ln \left( 1+e_{\binom{21}{32}%
}\right) $.

\section{Conclusions}

We have demonstrated that the dual group coordinates provide the natural
basis for the costructure of the twisted universal enveloping algebras. The $%
\mathfrak{g}\longrightarrow \mathfrak{g}^{\#}$-transition
applied to the twisted
algebra $U_{\mathcal{F}}\left( \mathfrak{g}\right) $ simplifies the problem of
finding solutions to the twist equations. In the forthcoming publication we
shall demonstrate how new classes of solutions are obtained in terms of $%
\mathfrak{g}^{\#}$-coordinates.

The dual group approach provides the new insight in the effect of the deformed
carrier spaces \cite{Kul-Lya},\cite{AKL}. In the weight system $\Gamma _{%
\mathcal{F}}^{\#}$ the weights $\lambda _{\perp }$ orthogonal to the initial
root $\lambda _{0}$ (of the full extended twist $\mathcal{F}$) cannot be
reached from the points of $\Gamma _{\mathfrak{a}}^{\#}$ by the shifts
corresponding to weights in $\Gamma _{\mathfrak{c}}^{\#}$. The reason is that
the system $\Gamma _{\mathfrak{c}}^{\#}$ is
located in the negative sector while
$\Gamma _{\perp \lambda _{0}}$ has the zero level, $\lambda _{\perp }\left(
h_{\lambda _{0}}\right) =0$. For the canonical extended twists (without the
Reshetikhin ''rotation'') this means that the coproducts $\Delta _{\mathcal{F%
}}(e_{\lambda _{\perp }})$ are primitive. This primitivity is realized only
in $\mathfrak{g}^{\#}$-coordinates.
Thus to find the additional twisting factors with the carrier
subalgebras in $V_{\perp } $ one must redefine $V_{\perp }$ in
terms of $\mathfrak{g}^{\#}$-basis.

The exponential basis used in this paper can not
be considered as universal. When the dual group is not solvable one must use
other dual group bases (for example, the matrix coordinates are to be used
for the parabolic dual group $\mathfrak{G}_{\mathcal{P}}^{\#}$ that contains
the simple subgroup in $\mathfrak{SL}(n)$).

\subsection{Acknowledgements}

The author is grateful to Prof. P.P.Kulish for stimulating discussions. The
work was supported by the Russian Foundation for Fundamental Research, grant
N 030100593.

\section{Appendix}

We can check that the limit $\lim_{\varepsilon ,\xi _{s}\longrightarrow
0;\quad \xi _{s}/\varepsilon =\zeta _{s}}U_{\mathcal{F}\left( \left\{ \xi
_{s}\right\} \right) }\left( \mathfrak{g}\left( \varepsilon \right) \right) $
exists and corresponds to the Hopf algebra $\mathrm{Fun}\left( \mathfrak{G}%
^{\#}\left( \left\{ \zeta _{s}\right\} \right) \right) $. It must be taken
into account that the twisting factors $\mathcal{F}_{q}$ can be deformed by
the previous twists \cite{Kul-Lya}. Their form can differ from the canonical
one, $\mathcal{F}_{\mathcal{J}}=\exp \left( h\otimes \sigma _{\mu }\right) $
for the Jordanian BTF and $\mathcal{F}_{\mathcal{E}}=\exp \left( e_{\lambda
}\otimes e_{\nu }f\left( \sigma \right) \right) $ for the extension. The
expressions $h\otimes \sigma _{\mu }$ and $e_{\lambda }\otimes e_{\nu
}f\left( \sigma \right) $ are the zero order terms in the expansion of the
deformed $\ln \mathcal{F}_{\mathcal{J}s}^{\mathrm{def}}$ and $\ln \mathcal{F}%
_{\mathcal{E}s}^{\mathrm{def}}$ \cite{L-BTF} with respect to the deformation
parameters of the preceding twists.

Let the parameters $\xi _{s}$ be proportional to $\varepsilon $
\begin{equation*}
\xi _{s}=\varepsilon \zeta _{s},
\end{equation*}
We shall assume that the logarithms $\ln \mathcal{F}_{q}$ of the twisting
factors $\mathcal{F}_{q}$ behave as follows
\begin{equation}
\ln \mathcal{F}_{q}=\frac{1}{\varepsilon }\Psi _{q}\left( \left\{ \widehat{%
e_{j}};\zeta _{s}\right\} \right) ,\qquad \Psi \in U\otimes U;\quad
s=1,\ldots ,l;q=1,\ldots ,p.  \label{lim-behav}
\end{equation}
It can be directly checked that in the type of deformation that we consider
here this condition always holds \cite{AKL}.

\begin{lemma}
In the Hopf algebra $U_{\mathcal{F}}\left( \mathfrak{g}\right) $ twisted by
the factorizable twist $\mathcal{F}=\mathcal{F}_{p}\mathcal{F}_{p-1}\ldots
\mathcal{F}_{1}$ the coproducts $\frac{\varepsilon }{\xi ^{\vartheta _{j}}}%
\Delta _{\mathcal{F}}\left( e_{j}\right) $ have the finite limit
\begin{equation*}
\lim_{\varepsilon ,\xi _{s}\longrightarrow 0;\quad \xi _{s}/\varepsilon
=\zeta _{s}}\left( \frac{\varepsilon }{\xi ^{\vartheta _{j}}}\Delta _{%
\mathcal{F}}\left( e_{j}\right) \right) =\Delta _{\mathcal{F},\left\{ \zeta
_{s}\right\} }^{\lim }\left( \widehat{e_{j}}\right) .
\end{equation*}
\end{lemma}

\begin{proof}
Any operator $\mathrm{ad}\left( \Psi _{q}\left( \left\{ \widehat{e_{j}}%
;\zeta _{s}\right\} \right) \right) $ with respect to $\varepsilon $ acts as
a multiplication by a polynomial with strictly positive powers (in
particular it multiplies by $\varepsilon $ any tensor belonging
to $ V_{\mathfrak{g}} \otimes U\left( \mathfrak{g}\right) $ or
$ U\left( \mathfrak{g}\right) \otimes V_{\mathfrak{g}} $ ).
In the limit the twisted
coproducts of the generators $\Delta _{\mathcal{F}}\left( e_{j}\right) $ are
defined by the action of the $\mathrm{BCH}\left\{ \Psi _{q}\left( \left\{
\widehat{e_{j}};\zeta _{s}\right\} \right) ;\varepsilon \right\} $ series
that has only positive powers of $\varepsilon $. Finally for any $e_{i}\in
\mathfrak{g}$ we get
\begin{equation*}
\begin{array}{l}
\lim_{\varepsilon ,\xi _{s}\longrightarrow 0;\quad \xi _{s}/\varepsilon
=\zeta _{s}}\left( \frac{\varepsilon }{\xi ^{\vartheta _{j}}}\Delta _{%
\mathcal{F}}\left( e_{j}\right) \right)  =\\
\rule{5mm}{0mm} =\lim_{\varepsilon ,\xi
_{s}\longrightarrow 0;\quad \xi _{s}/\varepsilon =\zeta
_{s}}\prod_{j=k}^{1}e^{\mathrm{ad}\left( \ln \mathcal{F}_{q}\right) }\circ
\frac{\varepsilon }{\xi ^{\vartheta _{m}}}\Delta ^{\left( 0\right) }\left(
e_{j}\right) = \\
\rule{5mm}{0mm}=\lim_{\varepsilon ,\xi _{s}\longrightarrow 0;
\quad \xi _{s}/\varepsilon
=\zeta _{s}}\left( e^{\frac{1}{\varepsilon }\mathrm{ad}\left( \mathrm{BCH}%
\left\{ \Psi _{q}\left( \left\{ \widehat{e_{j}};\zeta _{s}\right\} \right)
;\varepsilon \right\} \right) }\circ \frac{\varepsilon }{\xi ^{\vartheta
_{m}}}\Delta ^{\left( 0\right) }\left( e_{j}\right) \right) = \\
\rule{5mm}{0mm}=\Delta _{\mathcal{F},\left\{ \zeta _{s}\right\} }^{\lim }
\left( \widehat{e_{i}}\right) .
\end{array}
\end{equation*}
\end{proof}

\begin{corollary}
Notice that in the limit all terms arising due to rearrangement of the
monomials in the PBW basis fade away. Only the zero power term of $\mathrm{%
BCH}_{\mathcal{F}}\left\{ \Psi _{q}\left( \left\{ \widehat{e_{j}};\zeta
_{s}\right\} \right) ;\varepsilon \right\} $ will give the contribution to
the limit value $\Delta _{\mathcal{F},\left\{ \zeta _{s}\right\} }^{\lim
}\left( \widehat{e_{i}}\right) $,
\begin{equation}
\Delta _{\mathcal{F},\left\{ \zeta _{i}\right\} }^{\lim }\left( \widehat{%
e_{i}}\right) =\lim_{\varepsilon \longrightarrow 0}\left( e^{\frac{1}{%
\varepsilon }\mathrm{ad}\left( \sum_{q}\Psi _{q}\left( \left\{ \widehat{e_{j}%
};\zeta _{s}\right\} \right) \right) }\circ \left( \widehat{e_{i}}\otimes
1+1\otimes \widehat{e_{i}}\right) \right) .\square    \label{lim-group}
\end{equation}
\end{corollary}

Let the group $\mathfrak{G}^{H}\left( \left\{ \xi _{s}\right\} \right) $ be
defined by the commutative Hopf algebra $H$ with the costructure $\left\{
\Delta _{\mathcal{F},\left\{ \zeta _{s}\right\} }^{\lim }\right\} .$ This is
the second classical limit of $U_{\mathcal{F}}\left( \mathfrak{g}\right) $,
\begin{equation*}
H\left( \mathrm{Ab},\Delta _{\mathcal{F},\left\{ \zeta _{s}\right\} }^{\lim
}\right) \approx \mathrm{Fun}\left( \mathfrak{G}^{H}\left( \left\{ \zeta
_{s}\right\} \right) \right) .
\end{equation*}
It is sufficient to check that the corresponding Lie algebras are
equivalent.
The Lie coalgebra of the group $\mathfrak{G}^{H}\left( \left\{ \zeta
_{s}\right\} \right) $ is constructed in an ordinary way: we put all the
parameters proportional to $\xi $%
\begin{equation*}
\zeta _{s}=\xi _{s}\xi ,
\end{equation*}
and antisymmetrize the first power terms in the coproducts
\begin{equation*}
\delta \left( \widehat{e_{i}}\right) =\frac{d}{d\xi }\left( \Delta _{%
\mathcal{F},\left\{ \xi _{s}\xi \right\} }^{\lim }-\tau \circ \Delta _{%
\mathcal{F},\left\{ \xi _{s}\xi \right\} }^{\lim }\right) \left( \widehat{%
e_{i}}\right) \mid _{\xi =0}=
\end{equation*}
\begin{equation*}
=\frac{d}{d\xi }\mid_{\xi =0}
\lim_{\varepsilon \longrightarrow 0}\left( \left(
\begin{array}{c}
e^{\frac{1}{\varepsilon }\mathrm{ad}\left( \sum_{q}\Psi _{q}\left( \left\{
\widehat{e_{j}};\xi _{s}\xi \right\} \right) \right) }- \\
-e^{\frac{1}{\varepsilon }\mathrm{ad}\left( \sum_{q}\tau \circ \Psi
_{q}\left( \left\{ \widehat{e_{j}};\xi _{s}\xi \right\} \right) \right) }
\end{array}
\right) \circ \left( \widehat{e_{i}}\otimes 1+1\otimes \widehat{e_{i}}%
\right) \right) .
\end{equation*}
The functions $\Psi _{\left( \mathcal{F}\right) q}\left( \left\{ \widehat{%
e_{j}};\xi _{s}\xi \right\} \right) $ corresponding to the twisting factors $%
\mathcal{F}_{q}$ in $\mathcal{F}$ obey the rule:
\begin{equation*}
\frac{d}{d\xi }\mid_{\xi =0}\sum_{q}\left( \Psi _{\left( \mathcal{F}\right)
q}\left( \left\{ \widehat{e_{j}};\xi _{s}\xi \right\} \right) -\tau \circ
\Psi _{\left( \mathcal{F}\right) q}\left( \left\{ \widehat{e_{j}};\xi
_{s}\xi \right\} \right) \right) =r_{\mathcal{F}}\left( \left\{ \xi
_{s}\right\} \right) .
\end{equation*}
Thus the Lie coalgebra of the group $\mathfrak{G}^{H}\left( \left\{ \zeta
_{s}\right\} \right) $ described by $H\left( \mathrm{Ab},\Delta _{\mathcal{F}%
,\left\{ \zeta _{s}\right\} }^{\lim }\right) $ is fixed by the same
relations as those that
define $\mathfrak{g}^{\#}\left( \left\{ \xi _{s}\right\} \right) $ ,
\begin{equation*}
\delta \left( \widehat{e_{i}}\right) =\left[ r_{\mathcal{F}}\left( \left\{
\xi _{s}\right\} \right) ,\Delta ^{\left( 0\right) }\left( \widehat{e_{i}}%
\right) \right] .
\end{equation*}

The groups $\mathfrak{G}^{H}\left( \left\{ \xi _{s}\right\} \right) $
and $\mathfrak{%
G}^{\#}\left( \left\{ \xi _{s}\right\} \right) $ are simply connected and
have the same Lie algebra $\mathfrak{g}^{\#}\left( \left\{ \xi _{s}\right\}
\right) $.

\begin{remark}
It is essential to check the algebra $\mathfrak{g}^{H}\left( \left\{ \xi
_{s}\right\} \right) $. Varying the behaviour of the deformation parameters
we can involve the additional limit procedures and
obtain various contracted groups $\lim \mathfrak{G}^{\#}\left( \left\{ \xi
_{s}\right\} \right) $ corresponding to different boundaries of the
initial parametrized set. This is the uniform character of the parameters
with respect to $\xi $ that guarantees the isomorphism $\mathfrak{G}^{H}\left(
\left\{ \xi _{s}\right\} \right) \approx \mathfrak{G}^{\#}\left( \left\{ \xi
_{s}\right\} \right) $.
\end{remark}

\end{document}